\documentclass{amsart}
\usepackage{amsmath}
\usepackage{amsthm}
\usepackage{amssymb}
\usepackage{amsfonts}
\usepackage{enumerate}
\usepackage{multicol}
\usepackage[all]{xy}

\normalsize

\newtheoremstyle{theoremstyle}
  {10pt}      
  {5pt}       
  {\itshape}  
  {}          
  {\bfseries} 
  {:}         
  {.5em}      
  {}          

\newtheoremstyle{examplestyle}
  {10pt}      
  {5pt}       
  {}          
  {}          
  {\bfseries} 
  {:}         
  {.5em}      
  {}          

\theoremstyle{theoremstyle}
\newtheorem{theorem}{Theorem}[section]
\newtheorem{lemma}[theorem]{Lemma}
\newtheorem{proposition}[theorem]{Proposition}
\newtheorem{corollary}[theorem]{Corollary}
\theoremstyle{examplestyle}
\newtheorem{definition}[theorem]{Definition}
\newtheorem{remark}[theorem]{Remark}

\newcommand{\comment}[1]{}

\newcommand{\rays}{\Delta(1)}

\newcommand{\sh}[1]{\mathcal{#1}}
\newcommand{\msh}[1]{$\sh{#1}$}
\newcommand{\spec}[1]{\operatorname{spec}(#1)}
\newcommand{\tproj}[1]{\operatorname{tproj}(#1)}
\newcommand{\mtproj}[1]{$\tproj{#1}$}
\newcommand{\Spec}[1]{\operatorname{\mathbf{spec}}(#1)}
\newcommand{\Tproj}[1]{\operatorname{\mathbf{tproj}}(#1)}

\newcommand{\coord}[1]{k[{#1}_M]}
\newcommand{\mTCDiv}[1]{$\TCDiv{#1}$}
\newcommand{\TCDiv}[1]{\operatorname{TCDiv(#1)}}
\newcommand{\Supp}{\operatorname{SF}(N,\Delta)}
\newcommand{\mrSupp}{$\rSupp$}
\newcommand{\rSupp}{\operatorname{SF}(N,\Delta)_\mathbb{R}}
\newcommand{\pic}{\operatorname{Pic}}

\title{Toric Varieties as Spectra of Homogeneous Prime Ideals}
\author{Markus Perling}

\address{Universit\'e Grenoble, Institut Fourier, 100 rue des Maths, BP 74, 38402 St. Martin d'Heres, France}

\email{perling@mozart.ujf-grenoble.fr}


\date{June 2005}

\begin{document}

\begin{abstract}
We describe the construction of a class of toric varieties as spectra of
homogeneous prime ideals.
\end{abstract}
\maketitle

\section{Introduction}

Homogeneous coordinate rings for toric varieties, which are due to Cox \cite{Cox},
by now are a quite standard tool for the investigation of toric varieties and sheaves
over them (see \cite{Mustata3} for instance). It has been noticed earlier
\cite{ACampoHausenSchroeer}, that for a toric variety there may be many choices of
homogeneous coordinate rings and it is not a priori clear that the original
construction of Cox is the best choice. However, as shown in
\cite{ACampoHausenSchroeer}, Cox' homogeneous coordinate ring has the universal
property that for every other coordinate ring there exists a (suitably defined)
graded ring homomorphism into Cox' homogeneous coordinate ring. 

In this paper we consider the class of toric varieties which have enough invariant
Cartier divisors, first introduced by Kajiwara \cite{Kajiwara} (see also Definition
\ref{enoughdef}). We show that for any such toric variety $X$ there exists a
homogeneous coordinate ring $S = \bigoplus_{g \in G} S_g$, graded by a certain
abelian group $G$, such that $X$ can be identified with the set of homogeneous prime
ideals of $S$ minus a certain exceptional subset, along the lines of the
{\it proj}-construction in \cite{EGAII}.
The possibility of this construction is implicitely already contained in \cite{Cox}
and was remarked in \cite{Kajiwara} as a consequence of the constructions there
(see also \cite{BrennerSchroeer} for a related construction).
The importance of this construction from our point of view is that it allows us to
define quasicoherent sheaves over $X$ from $G$-graded $S$-modules by taking
homogeneous sections and without gluing.

We want to remark that we derived the notion of a toric variety with enough
invariant Cartier divisors without knowing about \cite{Kajiwara} and we hope
that our construction will provide a complementary point of view of Kajiwara's
constructcions.

\section{Basic Facts for Toric Varieties}

The objects we will study are {\it algebraic schemes over an algebraically closed
field $k$}, i.e. schemes $X$ such that there exists a morphism $X \longrightarrow
\Spec{k}$ of finite type. We will call such objects varieties and, if
not stated otherwise, we will assume that all varietes are reduced.
A toric variety $X$ is a normal variety which contains an algebraic torus
$T$ as an open dense subset such that the torus multiplication extends to
the action of the algebraic group $T$ on $X$. For background on toric varieties we refer
to standard literature such as \cite{Oda} and \cite{Fulton}. A toric variety
$X = X_\Delta$ is defined by a fan $\Delta$ contained in the real vector space
$N_\mathbb{R} \cong N
\otimes_{\mathbb{Z}} \mathbb{R}$ spanned over a lattice $N \cong \mathbb{Z}^d$.
We will always assume that a fan $\Delta$ is not contained in a proper subvector space of
$N_\mathbb{R}$. 
Let $M$ be the lattice dual to $N$ and let $\langle \, , \rangle : M \times N \rightarrow
\mathbb{Z}$ be the natural pairing
which extends to the scalar extensions $M_\mathbb{R} := M \otimes_\mathbb{Z} \mathbb{R}$
and $N_\mathbb{R}$. Elements of $M$ are denoted by $m$, $m'$, etc. if written additively,
and by $\chi(m)$, $\chi(m')$, etc. if written multiplicatively, i.e. $\chi(m + m') = \chi(m)
\chi(m')$. The lattice $M$ is the natural group of characters of the torus
$T = \operatorname{Hom}(M,k^*) \cong (k^*)^d$.

A cone $\sigma$ of the fan $\Delta$ is a convex rational polyhedral cone contained in
$N_\mathbb{R}$. For these cones the following standard notation will be used:

\begin{itemize}
\item cones are denoted by Greek letters $\rho$, $\sigma$, $\tau$, etc.,
the natural order among cones is denoted by $\tau < \sigma$;
\item $\Delta(1) := \{\sigma \in \Delta \mid \dim \sigma = 1\}$ is the set of
{\em rays} of $\Delta$;
\item $\vert \Delta \vert := \bigcup_{\sigma \in \Delta} \sigma$ is the {\it support} of the fan $\Delta$;
\item $n(\rho)$ denote the primitive lattice element corresponding to a ray $\rho$;
\item $\check{\sigma} := \{m \in M_\mathbb{R} \mid \langle m, n \rangle \geq 0
\text{ for all $n \in \sigma$}\}$
is the {\it dual} cone of $\sigma$;
\item $\sigma^\bot = \{m \in M_\mathbb{R} \mid m(n) = 0 \text{ for all } n \in \sigma \}$;
\item $\sigma_M := \check{\sigma} \cap M$ subsemigroup of $M$ associated to $\sigma$;
\item $U_\sigma$ denotes the affine toric variety associated to the cone $\sigma$,
and the semigroup ring $k[\sigma_M] = \bigoplus_{m \in \sigma_M} k \cdot \chi(m)$
is identified with the coordinate ring of $U_\sigma$;
\item for $\rho \in \rays$ denote $D_\rho$ the associated irreducible $T$-invariant
Weil divisor;
\item $\TCDiv{X_\Delta}$ denotes the $\mathbb{Z}$-module group freely
generated by the $T$-invariant Cartier divisors on $X_\Delta$.
\end{itemize}

Recall that dualizing a cone $\sigma$ induces an order-reversing correspondence
between the faces of $\sigma$ and the faces of $\check{\sigma}$. This correspondence
is given by $\tau \longleftrightarrow \check{\sigma} \cap \tau^\bot$.

By \cite{MiyakeOda}, Proposition 6.1 and \cite{Fulton}, \S 3.4,
there exists a short exact sequence:
\begin{equation}
\label{maindiagram}
0 \longrightarrow M \overset{\iota}{\longrightarrow} \TCDiv{X_\Delta} \longrightarrow
\pic{X_\Delta}
\longrightarrow 0
\end{equation}
i.e. the Picard group of $X_\Delta$ is generated by \mTCDiv{X_\Delta} modulo the set of
$T$-invariant principal divisors.

A {\it $\Delta$-linear support function} $h$ over $\Delta$ is a map from
$\vert \Delta \vert$ to $\mathbb{R}$ such that its restriction to some cone $\sigma$ is
a linear map and $h(\sigma \cap N) \subset \mathbb{Z}$. Such a map is determined by a
collection of characters $\{m_\sigma\}_{\sigma \in \Delta} \subset M$ which are
compatible in the way that $m_\sigma - m_{\sigma'} \in (\sigma \cap \sigma')^\bot$ for
all $\sigma, \sigma' \in \Delta$. We denote $\operatorname{SF}(N, \Delta)$ the group of
$\Delta$-linear support functions over $\Delta$, which is a
finitely generated, free $\mathbb{Z}$-module. The characters $\chi(m_\sigma)$ are
rational functions over $U_\sigma$ and by compatibility they glue, with respect to the
open cover $U_\sigma$, $\sigma \in \Delta$, to a Cartier divisor $D_h$ on $X_\Delta$.
The map $h \mapsto D_h$ defined this way is a bijection and yields an
identification of $\operatorname{SF}(N, \Delta)$ with \mTCDiv{X_\Delta}.
We denote $\rSupp := \operatorname{SF}(N, \Delta) \otimes_\mathbb{Z} \mathbb{R}$.
The elements of \mrSupp\ are piecewise $\mathbb{R}$-linear
functions $\vert \Delta \vert \longrightarrow \mathbb{R}$.

\section{Toric Varieties as Spectra of Homogeneous Prime Ideals}

In this section we are going to give a variation of Kajiwara's \cite{Kajiwara} method
to construct toric varieties which have the property to have
{\it enough invariant Cartier divisors} as geometric quotients of some
quasi-affine variety.
By applying $\operatorname{Hom}_\mathbb{Z}(\_, k^*)$ to the sequence (\ref{maindiagram}) we obtain:
\begin{equation}
\label{dualmaindiagram}
1 \longrightarrow G \longrightarrow \hat{T} \longrightarrow T \longrightarrow 1
\end{equation}
where $\hat{T} \cong \spec{k[\TCDiv{X}]}$ and $G \cong \spec{k[\TCDiv{X}]}$.
As \mTCDiv{X} is a free $\mathbb{Z}$-module $\hat{T}$ is a torus, 
whereas $G$ is a diagonalizable group scheme which in general must not necessarily be
reduced or irreducible. However, we will consider here only varieties where $\pic(X)$ is
free, which is true for the most cases of interest (see also remark \ref{gradedremark}).
Any
$\Delta$-linear support function can be considered as a character of the torus $\hat{T}$.
Consequently, we write $\chi(h)$ for some element $h \in \Supp$ when the group law is
considered as multiplication.
We proceed by constructing a strongly convex rational polyhedral cone inside the real vector
space $\rSupp$ and by explicitly
deriving conditions whether this cone is large enough such that it defines a quotient presentation
for $X$.

\begin{definition}
For a ray $\rho \in \rays$ and its primitive element $n_\rho$ define $H_\rho := \{h
\in \rSupp \mid h(n_\rho) \geq 0 \}$.
\end{definition}

\begin{proposition}\label{ratcone}
Let $\rho \in \rays$. Then
\begin{enumerate}[(i)]
\item\label{ratconeone} $H_\rho$ is a half space in \mrSupp.
\item\label{ratconetwo} For $\rho \neq \rho' \in \rays$: $H_\rho \neq H_{\rho'}$.
\item\label{ratconethree} The boundary of $H_\rho$ is a rational hyperplane.
\end{enumerate}
\end{proposition}

\begin{proof}
(\ref{ratconeone}): First note that $H_\rho$ is a nontrivial convex cone: there exists
at least one
$m \in M$ such that $m(n_\rho) > 0$, and if $h, h' \in H_\rho$, then $h + h' \in H_\rho$.
Also, if $h(n_\rho) > 0$, then $(-h)(n_\rho) < 0$, hence $-h \notin H_\rho$. We have
to show that $\partial H_\rho = \{h \in H_\rho \mid h(n_\rho) = 0\}$ is a hyperplane
in \mrSupp. To do this, we first investigate the subfan $\Delta'$ of $\Delta$ generated
by the set of maximal cones $\{\sigma_1, \dots, \sigma_n\}$ containing $\rho$.
The vector space
$\operatorname{SF}(N, \Delta')_\mathbb{R}$ is isomorphic to $\mathbb{R}^k$ for
some $k > 0$. We claim that the subspace of $\operatorname{SF}(N, \Delta')_\mathbb{R}$
consisting of the support functions vanishing at $n_\rho$ is isomorphic to
$\mathbb{R}^{k - 1}$.

{\it Proof of the claim:} we do induction on $n$. Let first $n = 1$. Then
$\operatorname{SF}(N, \Delta')_\mathbb{R}$ is isomorphic to $M_\mathbb{R} /
\sigma_1^\bot$. The vanishing of some $m \in M_\mathbb{R} / \sigma_1^\bot$
on $n_\rho$ is equivalent to the fact that $m$ is in the image of $\rho^\bot$
in $M_\mathbb{R} / \sigma_1^\bot$. Now $\sigma_1^\bot \subset \rho_1^\bot$,
and $\rho_1^\bot$ has codimension 1 in $M_\mathbb{R}$, hence it has codimension
1 in $M_\mathbb{R} / \sigma_1^\bot$.

Induction step: let $n > 1$. Consider the fan $\Delta'$ generated by the cones
$\sigma_1,$ $\dots,$ $\sigma_{n - 1}$. $\operatorname{SF}(N, \Delta')_\mathbb{R}$
is isomorphic to $\mathbb{R}^k$ for some $k > 0$. Now we add the cone
$\sigma_n$ to the fan $\Delta'$ and denote the resulting fan by $\Delta' \cup
\sigma_n$. By restriction we get a map: $\operatorname{res} : \operatorname{SF}
(N, \Delta' \cup \sigma_n)_\mathbb{R} \longrightarrow \operatorname{SF}(N,
\Delta')_\mathbb{R}$. Denote the subvector space of support functions of
$\operatorname{SF}(N, \Delta')$ which vanish at $n_\rho$ by $\operatorname{SF}(N,
\Delta')_0$. We get a fibre diagram:
\begin{equation*}
\xymatrix{
\operatorname{res}^{-1}(\operatorname{SF}(N, \Delta')_0 \cap \operatorname{res}(\operatorname{SF}(N,
\Delta' \cup \sigma_n))) \ar[rr]^-{\operatorname{res}} \ar@{^{(}->}[d] & & \operatorname{SF}(N, \Delta')_0
\ar@{^{(}->}[d] \\
\operatorname{SF}(N, \Delta' \cup \sigma_n) \ar[rr]^{\operatorname{res}} & &  \operatorname{SF}(N, \Delta')
}
\end{equation*}
We denote the preimage $\operatorname{res}^{-1}(\operatorname{SF}(N, \Delta')_0 \cap
\operatorname{res}(\operatorname{SF}(N, \Delta' \cup \sigma_n)))$ by $\operatorname{SF}
(N, \Delta' \cup \sigma_n)_0$. There exists at least one support function in
$\operatorname{SF}(N, \Delta' \cup \sigma_n)$ which does not vanish at $n_\rho$,
(i.e. an element of $M$ which is not orthogonal to $n_\rho$) and whose restriction
to $\Delta'$ does also not vanish at $n_\rho$. Therefore we know that the image of
$\operatorname{SF}(N, \Delta' \cup \sigma_n)$ does not completely lie in
$\operatorname{SF}(N, \Delta')_0$. Also note that the kernel of $\operatorname{res}$
is contained in $\operatorname{SF} (N, \Delta' \cup \sigma_n)_0$. Because of this and
because by induction
assumption $\operatorname{SF}(N, \Delta')_0$ has codimension 1 in $\operatorname{SF}
(N, \Delta')$, its preimage $\operatorname{SF} (N, \Delta' \cup \sigma_n)_0$ must also
have codimension 1 in $\operatorname{SF} (N, \Delta' \cup \sigma_n)$. This proves the claim.

Now consider $\Delta'$ as subfan of $\Delta$. Again we draw a fibre diagram:
\begin{equation*}
\xymatrix{
\operatorname{res}^{-1}(\operatorname{SF}(N, \Delta')_0 \cap \operatorname{res}(\operatorname{SF}(N, \Delta)))
\ar[rr]^-{\operatorname{res}} \ar@{^{(}->}[d] & & \operatorname{SF}(N, \Delta')_0 \ar@{^{(}->}[d] \\
\operatorname{SF}(N, \Delta) \ar[rr]^{\operatorname{res}} & &  \operatorname{SF}(N, \Delta')
}
\end{equation*}
with notation analogously as above. Again, choosing an $m \in M$ not vanishing at
$n_\rho$, we have that $\operatorname{SF}(N, \Delta)_0$ has codimension 1 in
$\operatorname{SF}(N, \Delta)$.

(\ref{ratconetwo}): For $\rho \neq \rho'$ we have $\rho^\bot \neq {\rho'}^\bot$.
Then choose
an $m \in \rho^\bot \setminus {\rho'}^\bot$. By definition $m \in \partial H_\rho$,
but $m \notin \partial H_{\rho'}$. By multiplying $m$ with an appropriate positive
or negative real number, we have still that $m \in \partial H_\rho$, but $m \notin
H_{\rho'}$, hence $m \in (H_\rho \setminus H_{\rho'}) \neq \emptyset$.

(\ref{ratconethree}): Note that we could have done the proof of (\ref{ratconeone})
without change in the scalar extension $\Supp \otimes_\mathbb{Z} \mathbb{Q}$
instead of $\Supp \otimes_\mathbb{Z} \mathbb{R}$ such that $\partial H_\rho$
restricted to $\Supp \otimes_\mathbb{Z} \mathbb{Q}$ is also a hypersurface in
this space.
\end{proof}

For the next two propositions we need the following notion:

\begin{definition}
\label{enoughdef}
Let $X$ be a toric variety associated to a fan $\Delta$. $X$ is said to have
{\it enough invariant Cartier divisors} if for every cone $\sigma \in \Delta$
there exists an effective $T$-invariant Cartier divisor whose support is
precisely the union of the $D_\rho$ with $\rho \in \Delta(1) \setminus \sigma(1)$.
\end{definition}

In the sequel we will always assume that $X$ has enough invariant Cartier divisors.

\begin{proposition}
\label{fulldim}
The intersection $\bigcap_{\rho \in \rays} H_\rho$ in \mrSupp\ is a strongly convex
polyhedral cone. It has full dimension in \mrSupp.
\end{proposition}

\begin{proof}
By Proposition \ref{ratcone}, (\ref{ratconethree}) we know that this intersection is a
rational polyhedron containing all support functions $h$ which satisfy $h(n_\rho) \geq
0$ for all $\rho \in \rays$. If $-h$ also lies in the intersection, it follows that
$h(n_\rho) = 0$ for all $\rho \in \rays$ and thus $h = 0$ as the $n_\rho$ span
$N_\mathbb{R}$. Hence the intersection contains no proper subspace of \mrSupp.
Because $X$ has enough effective divisors, it is easy to see that there exists a support
function $h$ which takes positive values values on each cone in $\rays$ and therefore is
contained in the interior of every half space $H_\rho$ and thus in the interior of
the intersection $\bigcap_{\rho \in \rays} H_\rho$.
So there exists an open neighbourhood of $h$ also contained in
$\bigcap_{\rho \in \rays} H_\rho$,
from which we conclude that the intersection has full dimension.
\end{proof}

\begin{definition}
Denote by $C$ the cone in $\rSupp\check{\ }$ dual to $\bigcap_{\rho \in \rays}
H_\rho$ and identify $\bigcap_{\rho \in \rays} H_\rho$ with $\check{C}$. For some
$\rho \in \rays$ we denote $l_\rho$ the primitive vector of the positive ray orthogonal
to $H_\rho$.
\end{definition}

\begin{proposition}
\label{onetoone}
$C$ is a strongly convex rational polyhedral cone. For every $\sigma \in \Delta$,
the primitive vectors $\{l_\rho \mid \rho \in \sigma(1)\}$ span a face of $C$.
\end{proposition}

\begin{proof}
As by Proposition \ref{fulldim}, $\check{C}$ is strongly convex and has full dimension,
the dual cone $C$ also is strongly convex and has full dimension.
Therefore there is a one-to-one
correspondence of the bounding hyperplanes $\partial H_\rho$ of $\check{C}$
and the $l_\rho$. By Proposition \ref{ratcone}, (\ref{ratconetwo}) we know that the
half spaces differ pairwise,
so their orthogonals do. Because $X$ has enough effective Cartier divisors,
we know that there exists for each $\rho \in \rays$ a support function $h$
vanishing on $\rho$ and being positive on all other rays. The support
function $-h$ is, as $h$, contained in $\partial H_\rho$, but negative on
all other rays, hence $\partial H_\rho$ is not contained in the union
$\bigcup_{\rho' \in \rays \setminus \{\rho\}} H_{\rho'}$. This implies that
the orthogonal of $\partial H_\rho$ is not contained in the interior of any
face of $C$ and thus the $l_\rho$ span the one-dimensional faces of $C$.
Now let $\sigma \in \Delta$ and consider the intersection $\Sigma := \check{C} \cap
\bigcap_{\rho \in \sigma(1)} \partial H_\rho$, which is clearly a face of
$\check{C}$. We have to show that
the dual face in $C$ is precisely spanned by the $l_\rho$, where
$\rho \in \sigma(1)$.
But as $X$ has enough invariant Cartier divisors, there exists an effective $T$-invariant
divisor whose support is precisely $\rays \setminus \sigma(1)$, so $\Sigma$ is
not contained in any $\partial H_\rho$, where $\rho \in \rays \setminus \sigma(1)$.
This implies that the face of $C$ dual to $\Sigma$ is precisely spanned by the
$l_\rho$ with $\rho \in \sigma(1)$.
\end{proof}

\begin{definition}
Let $\sigma \in \Delta$, then we denote $\hat{\sigma}$ the face of $C$ which
is generated by the $l_\rho$ with $\rho \in \sigma(1)$. We denote $\hat{\Delta}$ the
subfan of $C$ which is generated by the $\hat{\sigma}$, where $\sigma \in \Delta$.
We denote $\hat{X} := X_{\hat{\Delta}}$ the corresponding toric variety.
\end{definition}

Clearly, $\hat{X}$ is an open toric subvariety of the affine toric variety $U_C$
associated to the cone $C$ and therefore quasi-affine. Let $m \in M_\mathbb{R}$
which we can identify with a support function by the inclusion $\iota$
(see diagram (\ref{maindiagram}), respectively consider its scalar extension to
$\mathbb{R}$). Then for $n \in N_\mathbb{F}$, we have
$\langle m, n \rangle \geq 0$ iff $n$ is contained in some $\sigma$ such that $m \in
\check{\sigma}$. Via $\iota$, $\check{\sigma}$ then
becomes a subset of $\bigcap_{\rho \in \sigma(1)} H_\rho$, i.e. of the cone dual
to $\hat{\sigma}$.
Proposition \ref{onetoone} shows that $\hat{\Delta}$ has the property that
$\hat{\sigma} \cap \hat{\sigma}' \in \hat{\Delta}$ whenever $\sigma \cap \sigma'
in \Delta$. So, by dualizing, the map $\rSupp\check{\ }
\longrightarrow N_\mathbb{R}$ maps every cone $\hat{\sigma}$ to $\sigma$ and
thus induces a surjective map of fans $(\Supp\check{\ }, \hat{\Delta})
\longrightarrow
(N, \Delta)$. This map, respetively its induced toric surjection $\pi : \hat{X}
\longrightarrow X$, is a {\em quotient presentation} of $X$ in the sense
of \cite{ACampoHausenSchroeer}.

\begin{remark}
We want to point out that in the case where $\Delta$ is simplicial, our quotient
presentation coincides with the construction of Cox. In particular, $\hat{\Delta}$
then also is a simplicial fan. However, if $\Delta$ contains
a non-simplicial cone $\sigma$, then $\hat{\sigma} \in \hat{\Delta}$ also is
non-simplicial. This follows directly from the proof of \ref{onetoone}:
the intersection $\bigcap_{\rho \in \sigma(1)} \partial H_\rho$ contains all
support functions $h$ which vanish on $\sigma$. But by piecewise linearity,
$h\vert_\sigma = 0$ is equivalent to $h\vert_\rho = 0$ for $\rho \in I$ and
$I \subset \sigma(1)$ such that the $\{\rho \in I\}$ form a maximal linearly
independent subset. Then $\bigcap_{\rho \in \sigma(1)} \partial H_\rho
= \bigcap_{\rho \in I} \partial H_\rho$ for any such $I$. It follows that
$\dim \hat{\sigma} = \dim \sigma$ and thus $\hat{\sigma}$ is simplicial iff
$\sigma$ is. Moreover, it follows that the fan $\hat{\Delta}$ is combinatorially
equivalent to $\Delta$.
\end{remark}

\begin{definition}
We denote $S$ the coordinate ring of $U_C$. For every $\sigma \in \Delta$, we
denote $B_\sigma \subset S$ the ideal describing the complement of $U_{\hat{\sigma}}$
in $U_C$. We call $B := \sum_{\sigma \in \Delta}$ the {\em irrelevant ideal} and
$\mathbf{V}(B)$ the vanishing set of $B$ in $U_C$.
\end{definition}

Clearly, the ideals $B_\sigma$ are $\Supp$-graded and so is $B$. The variety
$\mathbf{V}(B)$ is the complement of the union of all $U_{\hat{\sigma}}$,
where $\sigma \in \Delta$. Note that for any affine toric variety $U_\sigma$ and
$\tau < \sigma$ some face, the reduced ideal of the closed subvariety $U_\sigma
\setminus U_\tau$ is generated over $k[\sigma_M]$ by the characters $m \in \sigma_M
\cap \tau^\bot$ which are contained in the relative interior of $\sigma_M \cap \tau^\bot$.
The semigroup $\tau_M$ is generated by $\sigma_M$ and $-m$ for any such $m$.

For $\sigma \in \Delta$, the ideal $B_\sigma$ is generated by integral support
functions whose associated $T$-Cartier divisor is effective and has support
precisely on $\rays \setminus \sigma(1)$. For such a support function $h$, we
denote $S_{\chi(h)} := (\{\chi(h)^n\}_{n \geq 0})^{-1} S$ the localization of $S$
at $\chi(h)$. It is easy to see that for any other such $h'$, the rings
$S_{\chi(h)}$ and $S_{\chi(h')}$ are naturally isomorphic. So once and for all,
for every $\sigma \in \Delta$
we choose a distinguished integral element $h_\sigma$ which is contained in the
relative interior of the cone in $\check{C}$ dual to $\hat{\sigma}$.

The closed set $U_C \setminus \hat{X}$ by
construction has codimension at least 2 in $U_C$, and thus $\Gamma(\hat{X}, \sh{O}_{\hat{X}}) = S$.
The dual action of $G$ (see diagram (\ref{dualmaindiagram})) on $S$ induces an isotypical
decomposition:
\begin{equation*}
S = \bigoplus_{\alpha \in \pic(X)} S_\alpha
\end{equation*}
which defines a $\pic(X)$-grading of $S$.
For $\sigma \in \Delta$, the ring $S_{\chi(h_\sigma)}$ is $\pic(X)$-graded as well.
The subring $S^G_{\chi(h_\sigma)}$ of $G$-invariants of $S_{\chi(h_\sigma)}$ then
precisely coincides with $S_{(\chi(h_\sigma))} = (S_{\chi(h_\sigma)})_0$,
its subring of degree zero elements.

\begin{lemma}
\label{nagspec}
Let $\sigma \in \Delta$. Then
$S_{\chi(h_\sigma)}^G \cong S_{(\chi(h_\sigma))} \cong \coord{\sigma}$.
\end{lemma}

\begin{proof}
By construction $h_\sigma$ is a support function on $\Delta$ which vanishes on
$\sigma$ and is positive on all other cones. The degree-zero subring
of $S_{\chi(h_\sigma)}$ consists of elements of the form $\frac{f}{(\chi(h_\sigma))^n}$
where $n \in \mathbb{N}$, $\deg_{\pic(X)}(f) = n \cdot
\deg_{\pic(X)}(\chi(h_\sigma))$. Now let $m \in \sigma_M$; then $m(n) \geq
0$ for all $n \in \sigma$.
 For a fixed such $m$
we can always choose $k \geq 0$ such that $(k \cdot h_\sigma + m)(n) \geq 0$ for all
$n \in \vert\Delta\vert$, hence $\chi(h_\sigma)^k \cdot \chi(m)$ is a regular function
on $V$. But then $\frac{\chi(m) \cdot \chi(h_\sigma)^k}{\chi(h_\sigma)^k} = \chi(m)
\in S_{\chi(h_\sigma)}$, hence the zero component of $S_{\chi(h_\sigma)}$
is spanned by the characters $\chi(m) \in \sigma_M \subset \check{C} \cap \Supp$.
\end{proof}

An immediate consequence is:

\begin{proposition}
The pair $(X, \pi)$ is a categorical quotient of the action of $G$ on $\hat{X}$.
\end{proposition}

Below, by the toric $\operatorname{proj}$ construction, we will show that $X$ indeed
is a geometric quotient of $\hat{X}$ by $G$.

\

\noindent {\bf The toric $\operatorname{proj}$ construction:}
Denote $\spec{S}$ the spectrum of prime ideals of $S$ and $\Spec{S} = (\spec{S},
\sh{O})$ the scheme whose underlying topological space is $\spec{S}$.
We consider $\hat{X}$ as the quasi-affine scheme $\Spec{S} \setminus \mathbf{V}(B) =
(\hat{X}, \sh{O}_{\hat{X}})$.
This scheme has the usual topology with closed sets $\mathbf{V}(\mathfrak{a})
= \{\mathfrak{p} \in \hat{X} \mid \mathfrak{a} \subset \mathfrak{p} \}$, where
$\mathfrak{a}$ runs over the ideals of $S$, and open sets
$D(f) = \spec{S} \setminus \mathbf{V}((f))$, where $f \in S$.

\begin{definition}
We set $\tproj{S} := \{\mathfrak{p} \in \hat{X} \text{ homogeneous} \}$
with the canonical inclusion $\tproj{S} \overset{i}{\hookrightarrow} \hat{X}$.
We endow \mtproj{S} with the induced topology and denote $\mathbf{V}_+(\mathfrak{a})
:= \mathbf{V}(\mathfrak{a}) \cap \tproj{S}$ and $D_+(f)
:= D(f) \cap \tproj{S}$ for $f \in S$.
\end{definition}

The following lemma is immediate:

\begin{lemma}
\label{homogeneousclosed}
Let $\mathfrak{a}$ be any ideal in $S$ and denote by $\mathfrak{a}'$ the homogeneous
ideal generated by the sets $\{\operatorname{pr}_\alpha(a) \mid a \in \mathfrak{a}\}_{\alpha
\in \pic(X)}$,
where the $\operatorname{pr}_\alpha : S \twoheadrightarrow S_\alpha$ are the
vector space projections. Then for any homogeneous ideal $\mathfrak{b}
\subset S$ we have: $\mathfrak{a} \subset \mathfrak{b}$ iff $\mathfrak{a}'
\subset \mathfrak{b}$. Thus $\mathbf{V}_+(\mathfrak{a})
=\mathbf{V}_+(\mathfrak{a}')$.
\end{lemma}

\begin{lemma}
\label{dint}
\begin{enumerate}[(i)]
\item $D_+(fg) = D_+(f) \cap D_+(g)$.
\item The set $\{D_+(f) \mid f \in B \text{ homogeneous}\}$ forms a
basis for the topology of \mtproj{S}.
\end{enumerate}
\end{lemma}

\begin{proof}
We only show the second claim.
Let $U = \tproj{S} \setminus \mathbf{V}_+(\mathfrak{a})$ be an open subset of
\mtproj{S}, where $\mathfrak{a}$ is any nonzero ideal in $S$, and $\mathfrak{a}$
can by Lemma \ref{homogeneousclosed} be chosen as homogeneous. We show that each
$\mathfrak{p} \in U$ is contained in some $D_+(h) \subset U$, $h \in B$ homogeneous.
Because $\mathfrak{p} \nsupseteq B$ there exists an $f \in B$ homogeneous such that
$\mathfrak{p} \in D_+(f)$. Also, because $\mathfrak{p} \nsupseteq \mathfrak{a}$,
there exists $g \in \mathfrak{a}$, such that $\mathfrak{p} \in D_+(g)$ and $D_+(g)
\; \cap \; \mathbf{V}_+(\mathfrak{a}) = \varnothing$. Set $h := fg$, then
$\mathfrak{p} \in D_+(h) \subset U$ by Lemma \ref{dint}.
\end{proof}

Note that in the case $X$ is affine, we have $G = \{1\}$, $B = S = k[\sigma_M]$,
and the basis of topology $\{D_+(f) \mid f \in B \text{ homogeneous}\}$ coincides with
the usual basis $\{D(f) \mid f \in S\}$.

We define another sheaf of rings $\sh{O}'$ on $\spec{S}$ as follows:
\begin{equation*}
U \mapsto \{ s \in \sh{O}(U) \mid s \text{ fulfills (') below} \}
=: \Gamma(U, \sh{O}')
\end{equation*}

(') $s$ is locally given by a fraction $\frac{f}{g}$ where $f$, $g$ are
homogeneous elements of $S$ s.t. $\operatorname{deg}(\frac{f}{g}) = 0$.

\

Using this sheaf we define the sheaf of rings on $\tproj{S}$ by $\sh{O}_{\tproj{S}}
:= i^{-1}\sh{O}'$.
\begin{definition}
We denote the ringed space $\Tproj{S} := (\tproj{S},\sh{O}_{\tproj{S}})$ the {\em toric
$\operatorname{proj}$} of $S$.
\end{definition}

The next step is to endow $\Tproj{S}$ with a scheme structure isomorphic to that
of $X$ by constructing a map of topological spaces $\Tproj{S} \longrightarrow X$
which induces an isomorphism of ringed spaces.
Let $\tau < \sigma$ and $h_\sigma$, $h_\tau$ two support functions vanishing
on $\sigma$ and $\tau$, respectively. Then we have two natural, commutative squares:
\begin{multicols}{2}[][5cm]
\xymatrix{
\tproj{S_{\chi(h_\tau)}} \ar@{^{(}->}[d] \ar[rr]^\cong & & D_+(\chi(h_\tau)) \ar@{^{(}->}[d] \\
\tproj{S_{\chi(h_\sigma)}} \ar[rr]^\cong & & D_+(\chi(h_\sigma))
}
\xymatrix{
\spec{S_{(\chi(h_\tau))}} \ar[rr]^\cong \ar@{^{(}->}[d] & & \spec{\coord{\tau}} \ar@{^{(}->}[d] \\
\spec{S_{(\chi(h_\sigma))}} \ar[rr]^\cong & & \spec{\coord{\sigma}}
}
\end{multicols}
In the left diagram $\tproj{S_{\chi(h_\sigma)}}$ denotes the homogeneous prime
ideal spectrum of $S_{\chi(h_\sigma)}$ (note that if $f \in B$ then $D_+(f)
\cap \mathbf{V}(B) = \emptyset$) and the maps are induced by the naturality
of localization; the second diagram is induced by the naturality of the isomorphism of
Lemma \ref{nagspec}.
Our aim now is to fill the gap between the two diagrams by horizontal arrows.

\begin{lemma}
\label{homunit}
Let $f \in S_{\chi(h_\sigma)}$ be some homogeneous non-unit. Then $f$ can be written as:
$f = u \cdot f'$ where $u$ is a unit in $S_{\chi(h_\sigma)}$ and $f'$ has degree zero.
\end{lemma}

\begin{proof}
First assume that $f = h$ for some support function $h \in \Supp$. Then: $h \vert_\sigma = m$
for some $m \in M$. The support function $h' := h - m$ vanishes on $\sigma$, hence $h' \in
\bigcap_{\rho \in \sigma(1)} \partial H_\rho$, hence it is a unit in $S_{\chi(h_\sigma)}$,
and $m$ has degree zero, so we write $h = h' + m$ or multiplicatively $\chi(h) = \chi(m) \cdot
\chi(h')$.
Now let $f$ of the form $\sum_i a_i \chi(h_i)$ for $a_i \in k$ and $h_i \in \Supp$.
Each of the $\chi(h_i)$ can be written as $\chi(m_i) \cdot \chi(h'_i)$ for $\chi(h'_i)$ units
and $\chi(m_i)$ of degree 0. Then write $f = \chi(h'_1) \cdot \sum_i a_i \chi(h'_i - h'_1)
\chi(m_i)$. The sum has degree zero.
\end{proof}

\begin{proposition}
\label{tprojisscheme}
Let $\chi(h_\sigma)$ a homogeneous generator of $B$. Then there is an isomorphism of
locally ringed spaces $\pi_\sigma : (D_+(\chi(h_\sigma)),$ $\sh{O}_{\tproj{S}}
\mid_{D_+(\chi(h_\sigma))})$ $\cong$ $\Spec{S_{(\chi(h_\sigma))}}$ which is given by
$\mathfrak{p} \mapsto (\mathfrak{p} \cdot S_{\chi(h_\sigma)}) \cap S_{(\chi(h_\sigma))}$.
\end{proposition}

\begin{proof}
Because for
$\mathfrak{p} \in D_+(\chi(h_\sigma))$, $\chi(h_\sigma)$ is not contained in $\mathfrak{p}$,
$\mathfrak{p} \cdot S_{\chi(h_\sigma)}$ is still a prime ideal in the localization. The
intersection of a prime ideal with a subring is also a prime ideal in this subring.
Hence, $(\mathfrak{p} \cdot S_{\chi(h_\sigma)}) \cap S_{(\chi(h_\sigma))}$ is prime
and $\pi_{\chi(h_\sigma)}$ is well-defined.
From Lemma \ref{homunit} we conclude that there exists a bijection between the sets
$\tproj{S_{\chi(h_\sigma)}}$ and $\spec{S_{(\chi(h_\sigma))}}$ and
$\pi_{\chi(f_\sigma)}(\mathfrak{p}) \subseteq
\pi_{\chi(f_\sigma)}(\mathfrak{p}')$ iff $\mathfrak{p} \subseteq \mathfrak{p}'$.
So there is a homeomorphism of topological spaces. We have still
to show that the stalks $\sh{O}'_\mathfrak{p}$ and $S_{{(\chi(h_\sigma))},
\pi_{\chi(h_\sigma)}(\mathfrak{p})}$ for
$\mathfrak{p} \in \tproj{S}$ are isomorphic. But the rings $S_{{(\chi(h_\sigma))},
\pi_{\chi(h_\sigma)}(\mathfrak{p})}$ are just elements $\frac{f}{g}$ where $f, g
\in S$ and $\deg_{\pic(X)}(f) = \deg_{\pic(X)}(g)$, and $g$ is not contained
in $\mathfrak{p}$.
\end{proof}

Let $\tau < \sigma$, then by naturality there is a commutative diagram:
\begin{equation*}
\xymatrix{
D_+(\chi(h_\tau)) \ar@{^{(}->}[d] \ar[rr]^{\pi_{\chi(h_\tau)}}_\cong & &
\spec{S_{(\chi(h_\tau))}} \ar@{^{(}->}[d] \\
D_+(\chi(h_\sigma)) \ar[rr]^{\pi_{\chi(h_\sigma)}}_\cong & & \spec{S_{(\chi(h_\sigma))}} \\
}
\end{equation*}

So we can conclude:

\begin{theorem}
The ringed space $\Tproj{S}$ has the structure of a scheme and is isomorphic to
the toric variety $X$.
\end{theorem}

\begin{corollary}
$X$ is a geometric quotient of $\hat{X}$ by $G$.
\end{corollary}

\begin{proof}
Let $\sigma \in \Delta$.
As seen in lemma \ref{homunit}, every homogeneous prime ideal in $U_{\hat{\sigma}}$ is
generated in degree zero, so there is a bijection between homogeneous prime ideals
and $\spec{k[\sigma_M]}$. In particular, to every geometric point $x \in U_\sigma$
corresponds a unique homogeneous ideal $\mathfrak{p}_x \in
\spec{S_{\chi(h_\sigma)}} \cap \tproj{S}$ which is maximal, and
thus $U_{\hat{\sigma}}$ is a fiber space
over $U_\sigma$. Now the result follows from \cite{Newstead}, Prop. 3.10.
\end{proof}

\begin{remark}\label{gradedremark}
The construction also works for the case where
$\pic(X)$ has a torsion subgroup. Then one has to modify slightly the notion
of homogeneous prime ideals and replaces them by so-called {\em $G$-prime} ideals.
i.e., an ideal $\mathfrak{p} \subset S$ is called $G$-prime if for any two
{\em homogeneous} elements $a, b \in S$ such that $ab \in \mathfrak{p}$ it follows
that either $a \in S$ or $b \in S$.
\end{remark}

Thanks to the $\Tproj{S}$ construction, we can define, quite in the fashion of
\cite{EGAII} for projective spaces, a quasicoherent sheaf
for a graded module without gluing:

\begin{definition}
Let $\mathfrak{p} \in \tproj{S}$ be a homogeneous prime ideal. Then we denote by
$T^{-1}_\mathfrak{p}$ the set of homogeneous elements of $S \setminus \mathfrak{p}$.
For any graded $S$-module $F$ the localization $T^{-1}_\mathfrak{p} F$ then is a graded
module as well. We denote by $F_{(\mathfrak{p})}$ the degree zero components of
$T^{-1}_\mathfrak{p} F$. For any open subset $U \subset \tproj{S}$ we define $\tilde{F}(U)$
to be the set of sections $s : U \longrightarrow \dot{\bigcup}_{\mathfrak{p} \in U}
F_{(\mathfrak{p})}$ such that $s(\mathfrak{p}) \in F_{(\mathfrak{p})}$ and for each
$\mathfrak{p} \in U$ there is a neighbourhood $V$ in $U$ and there are homogeneous
elements $m \in F$, $f \in S$ of same degree such that for each $\mathfrak{q} \in
V$ we have that $f \notin \mathfrak{q}$ and $s(\mathfrak{q}) = \frac{m}{f}$ in
$F_{(\mathfrak{q})}$. This defines a
presheaf on $\Tproj{S}$ and we denote by $\tilde{F}$ the associated sheaf.
\end{definition}

This sheafification operation defines an exact functor from the category of graded
$S$-modules to the category of quasicoherent sheaves over $X$. The sheaf $\tilde{F}$
is isomorphic to the sheaf obtained by the usual gluing procedure as in \cite{Kajiwara}.
In \cite{Kajiwara} it was shown that, if $F$ is finitely generated as an $S$-module,
then $\tilde{F}$ is a coherent sheaf.

\begin{proposition}
With notation as in the definition for any $\mathfrak{p} \in \tproj{S}$
there is an isomorphism $(\tilde{F})_{\mathfrak{p}} \cong F_{(\mathfrak{p})}$.
\end{proposition}

\begin{proof}
We define a map from $(\tilde{F})_{\mathfrak{p}}$ to $F_{(\mathfrak{p})}$ by mapping
$s \in (\tilde{F})_{\mathfrak{p}}$ to its value $s(\mathfrak{p})$. For any $\frac{m}{f}$
with $f \notin \mathfrak{p}$ there exists the neighbourhood $D_+(f)$ of $\mathfrak{p}$
such that $\frac{m}{f}$ exists in all $\mathfrak{q} \in D_+(f)$. Thus the map
is surjective. If there exist $s, t \in (\tilde{F})_{\mathfrak{p}}$ such that
$t(\mathfrak{p}) = s(\mathfrak{p})$, let's say $s(\mathfrak{p}) = \frac{m}{f}$ and
$t(\mathfrak{p}) = \frac{n}{g}$, then there exists $h \in T^{-1}_\mathfrak{p}$ such
that $h(gm - fn) = 0$. This is true for all $\mathfrak{q} \in D_+(f) \cap D_+(g) \cap
D_+(h)$. Thus $s = t$ in this neighbourhood, hence the germs of $s$ and $t$ in $\mathfrak{p}$
are identical. So the map is also injective, hence an isomorphism.
\end{proof}

In the other direction, for a quasicoherent sheaf \msh{F} on $X$ one can define a graded $S$-module
$\Gamma_*(X, \sh{F})$ $:= \bigoplus_{\alpha \in \pic{X}}$ $\Gamma(X, \sh{F}(\alpha))$, where
$\sh{F}(\alpha) := \sh{F} \otimes_{\sh{O}_X} \sh{O}_X(\alpha)$.
$\Gamma_*$ is a left-exact functor from the category of quasi-coherent sheaves over $X$ to the
category of graded $S$-modules. In \cite{Kajiwara} it is shown that there exists an isomorphism
$\Gamma_*(\sh{F})\tilde{\ } \cong \sh{F}$. In particular, every quasi-coherent $\sh{O}_X$--module
\msh{F} is of the form $\tilde{F}$ for some $\pic{X}$-graded $S$-module $F$. If \msh{F} is coherent,
then there exists a finitely generated $\pic{X}$-graded $S$-module $F$ with $\tilde{F} \cong \sh{F}$.

\end{document}